\documentclass[leqno, notitlepage]{amsart}
\usepackage{amscd}
\usepackage{amssymb}
\usepackage[all, poly]{xy}

\newcommand{\se}{\section}
\newcommand{\sus}{\subsection}
\newcommand{\sss}{\subsubsection}

\DeclareMathOperator{\Hom}{Hom}

\DeclareMathOperator{\supp}{supp}

\DeclareMathOperator{\Z}{\mathbb Z}

\DeclareMathOperator{\al}{\alpha}
\DeclareMathOperator{\la}{\lambda}

\DeclareMathOperator{\CLA}{\Lambda}
\DeclareMathOperator{\CQ}{Q}

\DeclareMathOperator{\GG}{\mathcal G}

\DeclareMathOperator{\kk}{{\bf k}}

\DeclareMathOperator{\IH}{{\bf IH}}

\DeclareMathOperator{\del}{\partial}
\DeclareMathOperator{\sm}{smooth}
\DeclareMathOperator{\Lie}{Lie}
\DeclareMathOperator{\Ad}{Ad}

\newcommand{\fg}{{\mathfrak g}}

\newcommand{\ca}{\check\alpha}

\newcommand{\cla}{\lambda}
\newcommand{\cmu}{\mu}

\newcommand{\com}{\check\omega}

\newcommand{\barr}{\overline}

\newcommand{\Fl}{{\mathfrak Fl}}

\newcommand{\OO}{{\mathcal O}}

\swapnumbers
\theoremstyle{plain} 
\newtheorem*{theorem}{Theorem}
\newtheorem*{corollary}{Corollary}

\newtheorem*{lemma}{Lemma}
\newtheorem*{proposition}{Proposition}

\newtheorem*{theoremA}{Theorem A}
\newtheorem*{theoremB}{Corollary B}

\numberwithin{equation}{subsubsection}

\numberwithin{enumi}{subsection}

\begin{document}

\title
{The Minimal degeneration singularities in the affine Grassmannians}

\author{ Anton Malkin }
\address{Department of Mathematics, MIT,
77 Massachusetts Avenue, Cambridge, MA 02139}

\email{malkin@math.mit.edu}

\author{ Viktor Ostrik }

\email{ostrik@math.mit.edu}

\author{ Maxim Vybornov }

\email{vybornov@math.mit.edu}


\begin{abstract} The minimal degeneration
singularities in the affine Grassmannians of simple simply-laced 
algebraic groups
are determined to be either Kleinian singularities of type A,
or closures of minimal orbits in nilpotent cones.
The singularities for non-simply-laced types are studied
by intersection cohomology and equivariant Chow group methods.
\end{abstract}

\maketitle

\se{Introduction}
\sus{} In the early 1980s H.~ Kraft and C.~ Procesi \cite{KP82}
classified the minimal degeneration singularities in the nilpotent
cones of classical Lie algebras stratified by the adjoint orbits.
We say that a pair of strata $\OO,\OO'$ in a variety
stratified by orbits of an algebraic group
is a \emph{minimal degeneration} if $\OO\subset \barr \OO'$,
and if $\OO\subseteq \barr\OO''\subseteq \barr\OO'$, 
then either $\OO=\OO''$ or
$\OO''=\OO'$.
The normal singularities of this sort turn out to be
smoothly equivalent to either Kleinian singularities of
type A and D, or the closures of minimal orbits. 
A similar problem for some flag varieties was more recently 
investigated by M.~Brion and
P.~ Polo \cite{BP}.

\sus{} In this paper we study the minimal degeneration
singularities for the $G[[z]]$-orbits in the 
affine Grassmannians. More precisely, 
let $G$ be a simple finite dimensional algebraic group 
over an algebraically closed field $\kk$ of characteristic $0$,
and let $\GG_G$ be the affine Grassmannian of $G$.
Using a Levi subgroup technique 
we obtain the following. 

\begin{theoremA} If $G$ is of simply-laced type,
then all singularities of the minimal degeneration singularities 
of $G[[z]]$-orbits in $\GG_G$ 
are either Kleinian singularities of type $A$ or
minimal singularities (closures of the minimal nilpotent orbits)
of types corresponding to Dynkin
subdiagrams of the Dynkin diagram of $G$.
\end{theoremA}

\sus{} 
For the non-simply-laced groups the situation is more
complicated. In addition to Kleinian and minimal singularities,
we obtain singularities which we call 
\emph{quasi-minimal}. These singularities are
studied by the methods of intersection cohomology \cite{L83},
and equivariant multiplicities, \cite{Ku, Br}.

\sus{} 
Our studies yield a completely new proof of the following result.

\begin{theoremB}
The smooth locus of the closure of a $G[[z]]$-orbit
in $\GG_G$ is precisely the $G[[z]]$-orbit itself.
\end{theoremB}

This result is due to S.~ Evens and I.~ Mirkovi\' c, it follows
immediately from \cite[Theorem 0.1(b)]{EM} which 
describes the characteristic cycle of the intersection 
cohomology sheaf on the  closure of a $G[[z]]$-orbit.

\sus{} The paper is organized as follows.
After introducing our notation and conventions
in the Section \ref{prelim}, we prove a ``Levi Lemma''
dealing with the action of the loop group of
a Levi subgroup of $G$ on $\GG_G$ 
in the Section \ref{LeviSection}.

In Section \ref{StemSection} we cite
a crucial result of J. Stembrige \cite{St} describing
the minimal degeneration of coweights
indexing the $G(O)$-orbits.

We then apply the Levi Lemma to describe
some minimal degeneration singularities
and prove Theorem A in the Section \ref{mindeglevi}.

Calculations of intersection cohomology dimensions
in the Section \ref{ic}
give us some information on the rational smoothness
of our singularities.

In Section \ref{rank2} we study the 
singularities appearing in the rank $2$ cases by methods of
\cite{Ku, Br}.

Finally, we deduce Corollary B in the section \ref{SmLocSec}.

\sus*{Acknowledgment}
We are grateful to A.~ Braverman,
M.~ Finkelberg, S.~ Kumar, I.~ Mirkovi\' c, and J.~ Starr 
for very useful discussions. The research of A.~M.~ and 
M.~V.~ was supported by the NSF Postdoctoral 
Research Fellowships. The research of  V.~O.~ was partially
supported by NSF grant DMS-0098830.

\se{Preliminaries}\label{prelim}
\sus{Notation}
\sss{}Let $G$ be a simple finite dimensional algebraic group
defined over an algebraically closed field $\kk$
of characteristic $0$.

Let us fix a Borel subgroup and a maximal torus 
$T\subset B\subset G$. 
We will
introduce the following notation and conventions:

\begin{enumerate}
\item{} $P\subset G$
is a paprabolic subgroup containing $B$.
\item{}$L_P=P/N_P$ is the Levi quotient
of $P$ by its nilpotent radical $N_P$.
\item{}$M=[L_P,L_P]$ is the commutant of  
$L_P$. $M$ is a connected semisimple group.
We will consider $M$ as a subgroup of $G$.
\item{} $I$ is the set of vertices of the Dynkin diagram
associated to $G$. The simple roots of the corresponding root system
are denoted by $\al_i$ for $i\in I$,
and the simple coroots by $\ca_i$ for $i\in I$.
\item{} $\CQ_G$ is the coroot lattice of $G$ and $\CQ^{\geq 0}_G$
is the nonnegative cone in $\CQ_G$.
\item{} $\CLA_G$ is the coweight lattice of $G$. 
The fundamental coweights are denoted
by $\com_i$ for $i\in I$. The partial order on the set of coweights:
$\cmu\geq\cla$ if and only if $\cmu-\cla\in\CQ^{\geq 0}_G$.
\end{enumerate}

\sss{} From now on we will assume that the parabolic 
$P$ is associated to a connected Dynkin subdiagram of the Dynkin
diagram of $G$. Then $M$ is a simple group of type
described by this Dynkin subdiagram. We will denote
the set of vertices of the subdiagram by $I_M$.
We need some more notation:
\begin{enumerate} 
\item{} $\CQ_M\subseteq Q_G$ is the coroot lattice of $M$, 
and $\CLA_M\subseteq \CLA_G$ is the coweight lattice of $M$.
\item{} Let $\la,\mu\in\CLA_G$. Then
$$
\cla=\sum_{i\in I}\cla_i\com_i \qquad \text{ and }
\qquad \cmu=\sum_{i\in I}\cmu_i\com_i .
$$
We denote:
 
\begin{equation}\label{lamum}
\cla_M=\sum_{i\in I_M}\cla_i\com_i \qquad \text{ and }
\qquad \cmu=\sum_{i\in I_M}\cmu_i\com_i .
\end{equation}
\end{enumerate}

\sus{Affine Grassmanians}

\sss{} Let $O=\kk[[z]]$ be the ring of formal power series
and $K=\kk((z))$ be its field of fractions. The 
affine Grassmannian
$\GG_G$ is the ind-scheme whose $\kk$-points are given by $G(K)/G(O)$.

We will recall here some facts about affine Grassmannians  
mostly borrowed from \cite{BL, BD, F, LS, NP, MVil}
where we refer the reader for many more details.

\sss{} The coweight $\cla\in\CLA_G$ may be 
considered as an element of $\GG_G$ 
via the identification $\CLA_G=\Hom(\kk^*, T)=T(K)/T(O)$.
It is well known
that the $G(O)$-orbits 
on $\GG_G$ are indexed by dominant coweights. The notation: 
$$
\GG_{\cla}=G(O)\cdot {\cla}\ \text{ for } \ 
\cla\in \CLA^{+}_G .
$$
It is well known that $\GG_{\cla}\subseteq\barr\GG_{\cmu}$
if and only if $\cla\leq\cmu$, and that
$
\dim \GG_{\cla}=\langle 2\rho,\cla \rangle
$,
where $2\rho$ is the sum of positive roots.

\sss{}  Let us consider the group ind-scheme
$G(\kk[z^{-1}])$, and let 
$L^{<0}G$ be the kernel
of the map $G(\kk[z^{-1}])\to G$ defined by $z^{-1}\mapsto 0$.
The following Lemma is well known.

\begin{lemma}\label{slice}
Consider $\cla\in\CLA^{+}_G$.
The orbit $L^{<0}G\cdot \cla$ is a transverse slice
to $\GG_{\cla}$ at the point $\cla$. In other words:
\begin{enumerate}
\item[(i)] $L^{<0}G\cdot \cla$ is locally closed in $\GG_G$.
\item[(ii)] The action map 
$G(O)\times (L^{<0}G\cdot \cla)\to \GG_G$
is an open embedding.
\item[(iii)] For any $\cla\leq \cmu\in\CLA^{+}_G$ 
$$
\dim (L^{<0}G\cdot \cla)\cap \barr \GG_{\cmu} =
\dim\GG_{\mu}-\dim\GG_{\cla} .
$$
\end{enumerate}
\end{lemma}

\begin{proof} (i) is clear.
(ii) According to \cite[Lemme 2.1]{NP} 
the multiplication morphism 
$
G(O)\times L^{<0} G\to G(K)
$
is an open embedding. Then the action
morphism
$
G(O)\times (L^{<0}G\cdot \cla)\to \GG_G
$
is an open embedding and therefore is \' etale.
This also proves (iii).
\end{proof}

\begin{lemma} Let $\cla \le \cmu$ be two dominant coweights. The scheme
$(L^{<0}G\cdot \cla)\cap \barr \GG_{\cmu}$
is reduced, irreducible and normal.
\end{lemma}

\begin{proof} Since 
$(L^{<0}G\cdot \cla)\cap \barr \GG_{\cmu}$ is a transverse slice 
to $\GG_{\cla}$, and $\barr \GG_{\cmu}$ is a normal variety
\cite{Ku2, Mat, Fa, Lit}, 
the variety $(L^{<0}G\cdot \cla)\cap \barr \GG_{\cmu}$
is also normal and irreducible. 
\end{proof}

\sss{} We will need one more lemma on the $\kk^*$-action.
The group $\kk^*$ acts on $\GG_G$ by ``loop rotations''
$$
z\mapsto sz\qquad \text{ for } s\in\kk^* .
$$
The following Lemma is lifted from \cite{MVyb03}.

\begin{lemma}\label{staraction} We have:
\begin{enumerate}
\item[(i)] $L^{<0}G\cdot \cla \cap \barr \GG_{\cmu}$
is $\kk^*$-invariant.
\item[(ii)] $L\in L^{<0}G\cdot\cla$
if and only if $\lim_{s\to\infty} L=\cla$.
\end{enumerate}
\end{lemma}

\sus{Kleinian and minimal singularities}

\sss{} We will think of a Kleinian singularity of type $A_p$
as an invariant theory quotient
$$
\kk^2/(\Z/(p+1))
$$
of the affine space $\kk^2$ by the cyclic group of order $p+1$,
cf. \cite[III.6.1]{S80}.

\sss{} Let $\fg=\Lie G$ be the Lie algebra of our group $G$. 
Let $\al_{\max}$ be the maximal root.
It is well known that $\ca_{\max}$ is the short dominant
coroot.
Let $v_{\max}\in \fg_{\max}$ be a highest weight vector.
The conjugacy class $\OO_{\min}=\Ad G\cdot v_{\max}$ 
is the \emph{minimal}
nilpotent orbit, and its closure 
$$
\barr\OO_{\min}(\fg)=\OO_{\min}\sqcup\{0\}
$$
is called the \emph{minimal singularity} of type $\fg$.
We will also index the minimal singularities by the Dynkin
diagrams corresponding to $\fg$ and by small letters
corresponding to capital letters indicating the type of $\fg$.
For example the minimal singularity of $sl_3$ may be referred to
as either $\barr\OO_{\min}(sl_3)$ or 
$\barr\OO_{\min}(\xymatrix{\circ\ar@{-}[r] &\circ})$  
or minimal singularity of type $a_2$. For many
more details we refer the reader to \cite{KP81, KP82}.

\sss{Minimal singularities and affine Grassmannians}

Let $\cla=0$ and 
$\ca_{\max}$ be the short dominant
coroot.
We are grateful to I.~ Mirkovi\' c for explaining
to us the following

\begin{lemma}\label{minimal} 
There exists an isomorphism of algebraic varieties.
$$ 
(L^{<0}G\cdot 0)\cap \barr \GG_{\ca_{\max}}
\simeq\barr\OO_{\min}(\fg) .
$$
\end{lemma}

\begin{proof} By the construction of \cite[Page 182]{BD}
the variety 
$$
\barr\GG_{\ca_{\max}}\simeq \{0\}\sqcup\OO_{\min}\sqcup P(\OO_{\min})
$$
is the disjoint union of $3$ $G$-orbits, where 
$P(\OO_{\min})\subset P(\fg)$ the the projectivization of $\OO_{\min}$.
Now the variety $(L^{<0}G\cdot 0)\cap \barr \GG_{\ca_{\max}}$
is $G$-invariant, and therefore a union of $G$-orbits.
Since $P(\OO_{\min})$ is $\kk^*$-invariant, closed, and 
separated from $0$, by 
Lemma \ref{staraction} we have
$$
(L^{<0}G\cdot 0)\cap 
\barr \GG_{\ca_{\max}}\simeq \{0\}\sqcup\OO_{\min}=
\barr\OO_{\min} .
$$
\end{proof}

\se{The Levi Lemma}\label{LeviSection}

\sus{} Let us consider two dominant coweights
$\cla\leq\cmu \in \CLA^{+}_G$ such that their
difference $\cmu-\cla\in \CLA_M$ is in the
coroot lattice of $M$
generated by $\ca_i$ for $i\in I_M$ and therefore
\begin{equation}\label{diff}
\cmu-\cla=\cmu_M-\cla_M\in \CQ_M ,
\end{equation}
where $\cla_M, \ \cmu_M$ are defined in (\ref{lamum})

\sss{} Since we consider $M$ to be a subgroup of $G$,
the group $M(K)$ acts on the affine Grassmannian $\GG_G$.

\begin{lemma}\label{mm} 
There exists a natural ind-scheme isomorphism
$$
M(K)\cdot\cla \simeq \GG_M
$$
given by the map $m\cdot\cla\mapsto m\cdot\cla_M$
for $m\in M(K)$. Moreover, this isomorphism restricts
to the isomorphisms of varieties
$$
\GG_G\supset (L^{<0}M\cdot \cla)\cap \barr {M(O)\cdot \mu}
\simeq (L^{<0}M\cdot \cla_M)\cap \barr {M(O)\cdot \cmu_M} 
\subset \GG_M .
$$
\end{lemma}

\begin{proof} Indeed, 
$
M(K)\cdot\cla=M(K)/M(K)_{\cla}
$
where $M(K)_{\cla}$ is the stabilizer of $\cla$ in $M(K)$.
Now 
\begin{equation}\nonumber
\begin{split}
M(K)_{\cla}& =M(K)\cap z^{\cla} G(O) z^{-\cla} \\
& = z^{\cla}(z^{-\cla}M(K)z^{\cla}\cap G(O))z^{-\cla} \\
& =z^{\cla}(M(K)\cap G(O))z^{-\cla} \\
& =z^{\cla_M} M(O) z^{-\cla_M} .
\end{split}
\end{equation}
\end{proof}

\sss{} We will need another 

\begin{lemma}\label{mg} 
Consider $\cla\leq\cmu \in \CLA^{+}_G$
with $\cmu-\cla\in \CQ_M$. Then 
$$
L^{<0}M\cdot \cla\cap \barr {M(O)\cdot \cmu}
=L^{<0}G\cdot \cla\cap \barr {G(O)\cdot \cmu}
\subset \GG_G .
$$
\end{lemma}

\begin{proof} Denote 
$$
Y_M=L^{<0}M\cdot \cla\cap \barr {M(O)\cdot \mu}\ 
\text{ and }\
Y_G=L^{<0}G\cdot \cla\cap \barr {G(O)\cdot \mu} . 
$$
According to the lemma \ref{staraction}
both $Y_M$ and $Y_G$ are $\kk^*$-invariant
and therefore their closures $\barr Y_M$ and 
$\barr Y_G$ are also $\kk^*$-invariant.
Now due to the equation (\ref{diff}) and the Lemma \ref{mm}
$\dim Y_M=\dim Y_G$. Thus $\barr Y_M=\barr Y_G$. 

By construction, $Y_M\subseteq Y_G$. Now 
$\del Y_M=\barr Y_M-Y_M$ is closed since 
$Y_M$ is locally closed.

Now, consider a point $p\in Y_G-Y_M\subset \del Y_M$. According
to the Lemma \ref{staraction}, $\kk^*$-action retracts $p$ 
to $\cla$ but this is impossible since $p$ is in 
a closed and $\kk^*$-invariant set $\del Y_M$ which
does not contain $\cla$. This contradiction shows that
$Y_G-Y_M\cap \del Y_M=\emptyset$. 

Thus we have a bijective map of normal varieties $Y_M\to Y_G$
and thus an isomorphism.
\end{proof}

\begin{corollary}\label{levi}{\rm {\bf [Levi Lemma]}} 
There exists an isomorphism of algebraic varieties
$$
\GG_M\supset L^{<0}M\cdot \cla_M\cap \barr {M(O)\cdot \cmu_M}
=L^{<0}G\cdot \cla\cap \barr {G(O)\cdot \cmu}\subset \GG_G .
$$
\end{corollary}

\begin{proof} Follows from Lemma \ref{mm} and Lemma \ref{mg}.
\end{proof}

\se{Minimal degenerations of coweights}\label{StemSection}

\sus{} We will say that a pair $\cla, \cmu\in \CLA^{+}_G$ 
of dominant coweights
is a minimal degeneration if
\begin{enumerate}
\item[(i)] $\cla< \cmu$,
\item[(ii)] if $\check\nu\in\CLA^{+}_G$ is such that
$\cla\leq \check\nu\leq\cmu$ then either $\check\nu=\cla$
or $\check\nu=\cmu$.
\end{enumerate}

The pair $\cla,\cmu$ which is a minimal degeneration will be denoted
as $\cmu\leadsto\cla$.

\sus{} All minimal degenerations $\cmu \leadsto \cla$ are classified by
J.~Stembridge. We will say that the support $\supp(\cmu -\cla)$ of 
$\cmu -\cla$ is the Dynkin subdiagram involving all simple 
coroots appearing 
in the decomposition of $\cmu-\cla$. 
It is obvious that for a minimal degeneration 
$\supp(\cmu -\cla)$ is connected. Here is the Stembridge's list, 
see \cite[Theorem 2.8]{St}:

\begin{theorem}\label{Stem}\cite{St}
The pair $\cmu \leadsto \cla$ is a minimal degeneration
if and only if one of the following holds
\begin{enumerate}
\item \label{St1} $\cmu -\cla$ is a simple coroot $\al_i$, $i\in I$.
\item \label{St2} 
$\cmu -\cla$ is the short dominant coroot of $\supp(\cmu-\cla)$ and
$\cla =0$ on $\supp(\cmu -\cla)$.
\item \label{sp}
$\cmu -\cla$ is the short dominant coroot of $\supp(\mu-\lambda)$;
$\supp(\cmu-\cla)$ is of type $C_n$ and $\cla$ on 
$\supp(\cmu-\cla)$ is

$$
\xymatrix{
0 \ar@{-}[r] &
0 \ar@{-}[r] &
\cdots \ar@{-}[r] &
0 \ar@{<=}[r] &
1 } .
$$ 
\vskip .1in
\item \label{Stag2} $\supp(\cmu -\cla)$ is of type $G_2$, and 
$\cla =(\xymatrix{2\ar@3{>}[r] &0})$,
$\cmu =(\xymatrix{1\ar@3{>}[r] &1})$ on $\supp(\cmu -\cla)$.
\item \label{Stcg2} $\supp(\cmu -\cla)$ is of type $G_2$, and 
$\cla =(\xymatrix{1\ar@3{>}[r] &0})$,
$\cmu =(\xymatrix{0\ar@3{>}[r] &1})$ on $\supp(\cmu -\cla)$.
\end{enumerate}
\end{theorem}

\se{Minimal degenerations of $G(O)-$orbits 
and Levi subgroups}\label{mindeglevi}

\sus{The $PGL_2$ case}

\sss{} In the $PGL_2$ case the pair $\cla,\cmu \in\CLA^{+}_G$ 
is a minimal degeneration if and only if $\cmu=(p+2)\ \com$
and $\cla=p\ \com$, where $\com$ is the fundamental coweight and
$p\geq 0$. In this case the Main Theorem of \cite{MVyb}
implies the following

\begin{lemma}\label{sl2} Let $G=PGL_2$, $\cla=p\ \com$ and 
$\cmu=(p+2)\ \com$.
Then the singularity of the Schubert variety $\barr \GG_{\cmu}$
along $\GG_{\cla}$ is a Kleinian singularity of type $A_{p+1}$.
More precisely,
$$
(L^{<0}G\cdot \cla)\cap \barr \GG_{\cmu}\simeq 
\kk^2/(\Z/(p+2)) .
$$
\end{lemma}

\begin{proof} It is shown in \cite{MVyb}
that $(L^{<0}G\cdot \cla)\cap \barr \GG_{\cmu}$ is isomorphic
to a transverse slice to the subregular orbit in 
the nilpotent cone ${\mathcal N}\subset sl_{p+2}$.
The statement follows by a celebrated result
of Brieskorn and Slodowy \cite{S80}. For $p=0$ 
the lemma follows already from \cite{L81}.
\end{proof}

\sus{The proof of Theorem A}
\begin{theorem}\label{LeviTheorem} 
Let $\cla, \cmu\in\CLA^{+}_G$ be two
dominant coweights, let $\cla=\sum_{i\in I}\cla_i\com_i$,
and let $\cmu \leadsto \cla$ be 
a minimal degeneration. Then 
\begin{enumerate}
\item \label{case1} If $\cmu -\cla=\ca_i$ is a simple coroot for $i\in I$,
then
$$
(L^{<0}G\cdot \cla)\cap \barr \GG_{\cmu}\simeq \kk^2/(\Z/(\cla_i+2)).
$$
In other words, we have a Kleinian singularity of type $A_{\cla_i+1}$.
\item \label{case2} If $\cmu -\cla$ is the short dominant coroot of 
of the root system of type $\supp(\cmu-\cla)$ and
$\cla =0$ on $\supp(\cmu -\cla)$, then
$$
(L^{<0}G\cdot \cla)\cap \barr \GG_{\cmu}
\simeq\barr\OO_{\min}(\supp(\cmu-\cla)), 
$$
i.e., we have the minimal singularity of type $\supp(\cmu-\cla)$.
\end{enumerate}
\end{theorem}

\begin{proof} In the case (\ref{case1}) the Theorem follows
from the Levi Lemma (Corollary \ref{levi}) and
Lemma \ref{sl2}. In the case (\ref{case2}) the Theorem follows
from the Levi Lemma (Corollary \ref{levi}) and
Lemma \ref{minimal}. In the simply-laced case Theorem A follows.
If $G$ is of type $A$ the theorem follows already from \cite{MVyb}.
\end{proof} 

\se{Intersection cohomology calculations}\label{ic}

\sus{Notation}
\sss{Definition} We will say that a variety $X$ is rationally 
smooth at the point $x$ if the stalk of the intersection cohomology 
complex is $1$-dimensional in degree $(-\dim X)$, i.e. 
$\IH_x(X)=\kk[\dim X]$. 
We will say that $X$ is rationally smooth if it is rationally smooth
at every point, i.e $\IH(X)$ is the constant sheaf shifted by $\dim X$.

\sss{} By Lemma \ref{slice}, if $\cmu\leadsto\cla$ is a minimal 
degeneration,
the variety $(L^{<0}G\cdot \cla)\cap \barr \GG_{\cmu}$
is smooth (and rationally smooth) at every point but $\cla$.
Let us denote:
$$ 
m_{\cla}(\cmu,q)=\sum_i\dim 
\IH^{i-\dim(L^{<0}G\cdot \cla)\cap \barr \GG_{\cmu}}_{\cla}
(\barr \GG_{\cmu})\cdot q^{i} ,
$$
where $\IH^i(\barr \GG_{\cmu})$ is the $i$-th 
cohomology sheaf of the intersection
cohomology complex of $\GG_{\cmu}$ and  
$\IH^i_{\cla}(\barr \GG_{\cmu})$ is 
the stalk of $\IH^i(\barr \GG_{\cmu})$ at the point 
$\cla \in \barr \GG_{\cmu}$.
We will also consider the Euler characteristic
$$
\chi_{\cla}(\cmu)=m_{\cla}(\cmu,1).
$$
Clearly, $(L^{<0}G\cdot \cla)\cap \barr \GG_{\cmu}$ is 
rationally smooth
if and only if $m_{\cla}(\cmu)=\chi_{\cla}(\cmu)=1$.

\sss{} 
Consider $\cla$ and $\cmu$ as the dominant weights for the Langlands dual 
group $G^\vee$ and let $m_{\cla}(\cmu)$ denote the multiplicity of the weight
$\cla$ in the simple $G^\vee-$module with the highest weight $\cmu$. 
According to \cite{L83} one has
\begin{equation}\label{eulermult}
\chi_{\cla}(\cmu)=m_{\cla}(\cmu).
\end{equation}

\sus{Rational smoothness of minimal degenerations}

\sss{} Let us study the rational smoothness of the variety
$(L^{<0}G\cdot \cla)\cap \barr \GG_{\cmu}$ 
in the four cases of the Theorem \ref{Stem}
using the formula (\ref{eulermult}).

\begin{proposition}\label{ICProposition} 
Let $\cmu\leadsto\cla$ be a minimal 
degeneration and
\begin{enumerate}
\item  let $\cmu\leadsto\cla$ be as in \ref{St1}. Then
$$
m_{\cla}(\cmu,q)=1
$$
and the variety $(L^{<0}G\cdot \cla)\cap\barr\GG_{\cmu}$ 
is rationally smooth.
\item  let $\cmu\leadsto\cla$ be as in \ref{St2} with $\cla=0$.
Then 
\begin{equation}\nonumber
m_{\cla}(\cmu,q)=
\begin{cases}
\sum_{i=1}^n q^{e_i-1}, & \supp(\cmu-\cla) \text{ of type } ADE ,\\
\sum_{i=0}^{n-2} q^{2i},& \supp(\cmu-\cla) \text{ of type } B_n ,\\
1, & \supp(\cmu-\cla) \text{ of type } C_n ,\\
1+q^4, &  \supp(\cmu-\cla) \text{ of type } F_4 ,\\ 
1, & \supp(\cmu-\cla) \text{ of type } G_2 ,
\end{cases}
\end{equation}
where $e_i$ $1\leq i\leq n$ are the exponents for types ADE.
The Euler characteristic
$\chi_{\cla}(\cmu)=n$, $n-1$, $1$, $n$, $6$, 
$7$, $8$ ,$2$, $1$ for  $\supp(\cmu-\cla)$ of type 
$A_n$, $B_n$, $C_n$, $D_n$, $E_6$, $E_7$, $E_8$, $F_4$, $G_2$ 
respectively. Thus the variety 
$(L^{<0}G\cdot \cla)\cap\barr\GG_{\cmu}$ is rationally smooth
in the cases $C_n$ and $G_2$.
\item  \label{rationalacl} 
let $\cmu\leadsto\cla$ be as in \ref{sp}.
Then
$$
m_{\cla}(\cmu,q)=\sum_{i=0}^{n-1} q^i
$$
and $\chi_{\cla}(\cmu)=n$, and the variety 
$(L^{<0}G\cdot \cla)\cap\barr\GG_{\cmu}$ 
is not rationally smooth.
\item \label{rationalag2}  
let $G$ be of type $G_2$ and 
$\cla =(\xymatrix{2\ar@3{>}[r] &0})$, and
$\cmu =(\xymatrix{1\ar@3{>}[r] &1})$. 
Then
$$
m_{\cla}(\cmu,q)=1+q
$$
and $\chi_{\cla}(\cmu)=2$
and the variety 
$(L^{<0}G\cdot \cla)\cap\barr\GG_{\cmu}$
is not rationally smooth. 
\item  \label{rationalbg2} 
let $G$ be of type $G_2$ and 
$\cla =(\xymatrix{1\ar@3{>}[r] &0})$,
$\cmu =(\xymatrix{0\ar@3{>}[r] &1})$.
Then $m_{\cla}(\cmu,q)=\chi_{\cla}(\cmu)=1$
and the variety 
$(L^{<0}G\cdot \cla)\cap\barr\GG_{\cmu}$
is rationally smooth. 
\end{enumerate}
\end{proposition}

\begin{proof} 
First of all we apply the Levi Lemma
(Corollary \ref{levi}) to reduce
the statements to the Levi subgroup associated to
$\supp(\cmu-\cla)$.
Then the formulas are obtained by the application
of the results of \cite{L83} such as
the formula (\ref{eulermult}),
and the direct calculations.
\end{proof}

\sus{The quasi-minimal singularities}

\sss{Type $ac_n$} Let us look again at the case \ref{rationalacl}. 
The singular variety $(L^{<0}G\cdot \cla)\cap \barr \GG_{\cmu}$
has dimension $2n$, the same as the minimal singularity of type $a_n$,
and moreover, it follows from the Lemma that the polynomials
$m_{\cla}(\cmu,q)$ coincide in our case and the case of the
minimal singularity of type $a_n$:
$$
m_{\cla}(\cmu,q)=m^{a_n}_{\cla}(\cmu, q)=\sum_{i=0}^{n-1} q^i .
$$
We will call this
singularity arising from the affine Grassmannian of type $C_n$
the \emph{quasi-minimal} singularity of type $ac_n$.
Notice that according to Lemma \ref{sl2} 
the singularity $ac_1$ is the Kleinian singularity of type $A_2$.

\sss{Type $ag_2$} Let us look again at the case \ref{rationalag2}.
The singular variety $(L^{<0}G\cdot \cla)\cap \barr \GG_{\cmu}$
is $4$-dimensional as is the minimal singularity of type $a_2$
and moreover, it follows from the Lemma that the polynomials
$m_{\cla}(\cmu,q)$ coincide in our case and the case of the
minimal singularity of type $a_2$:
$$
m_{\cla}(\cmu,q)=m^{a_2}_{\cla}(\cmu, q)=1+q .
$$
We will call the
singularity arising from the affine Grassmannian  
the \emph{quasi-minimal} singularity of type $ag_2$.

\sss{Type $cg_2$} Let us look again at the case \ref{rationalbg2}.
The singular variety $(L^{<0}G\cdot \cla)\cap \barr \GG_{\cmu}$
is $4$-dimensional as is the minimal singularity of type $c_2$
and moreover it follows from the Lemma that  
both our variety and the minimal singularity of type $c_2$ 
are rationally smooth
(but not smooth, see Section \ref{rank2}).
We will call the
singularity arising from the affine Grassmannian  
the \emph{quasi-minimal} singularity of type $cg_2$.

\se{Equivariant multiplicities in rank 2}\label{rank2}

\sus{} In this section we will study 
the rank $2$ situation
in more detail.
Namely, for a minimal degeneration $\cmu \leadsto \cla$ the variety 
$(L^{<0}G\cdot \cla)\cap \barr \GG_{\cmu}$ is invariant under the action
of torus $T\times \kk^*$ and the point $\cla$ is fixed by this action. 
Here $T$ is the maximal torus of $G$ and $\kk^*$
acts by loop rotations. Thus
the {\em equivariant multiplicity}  
$e_{\cla}(\cmu)$ of the variety $(L^{<0}G\cdot \cla)\cap \barr \GG_{\cmu}$ at 
the point $\cla$ 
(i.e. the localization of the fundamental
class in the $(T\times \kk^*)-$equivariant Chow group at the only 
fixed point
$\cla$)
is defined, see e.g. \cite{Br}. Recall that $e_{\cla}(\cmu)$ 
is a rational function on the Lie algebra of the torus $T\times \kk^*$. In 
this section we calculate $e_{\cla}(\cmu)$ in all rank 2 cases when the 
codimension of degeneration is $>2$.  

Our calculation is performed as follows. Let $\Fl$ denote the affine flag 
variety of $G$ and let $\pi :\Fl \to \GG$ be the canonical projection. The map
$\pi$ is smooth with all of its fibers isomorphic to the finite dimensional 
flag variety $G/B$. Let $\GG_{\cla}\subset \barr \GG_{\cmu}$ be a minimal
degeneration. Then the singularity $\pi^{-1}(\GG_{\cla})\subset 
\pi^{-1}(\barr \GG_{\cmu})$ is smoothly equivalent to the singularity
$\GG_{\cla}\subset \barr \GG_{\cmu}$. Recall that the Schubert varieties in
$\Fl$ are labeled by the elements of the (extended) affine Weyl group 
$W_{aff}$. Let $y\in W_{aff}$ (respectively $w\in W_a$) label the open 
Schubert variety $X_y$ in $\GG_{\cla}$ (respectively $\GG_{\cmu}$). A formula 
for calculation of the equivariant multiplicity $e_yX(w)$ is given in 
\cite{Br} and one deduces easily from this a formula for equivariant 
multiplicity $e_{\cla}(\cmu)$ of the transversal slice, see 
\cite[page 27]{Br}
We perform our calculations using this formula and a computer.

\sus{Kumar's Criterion} In \cite{Ku}
S.~Kumar gave a general criterion for smoothness of Schubert varieties of
a general Kac-Moody group in terms of equivariant multiplicities. 
In particular, Kumar's criterion implies that if the variety
$\GG_{\cla}\subset \barr \GG_{\cmu}$ is smooth at $\cla$ then 
$$e_{\cla}(\cmu)=\prod_{\alpha \in S}\alpha^{-1},$$
where $S$ is a certain finite subset of the set of roots. We will see that
Kumar's criterion and our calculations imply that in all cases considered
in this section the varieties $(L^{<0}G\cdot \cla)\cap \barr \GG_{\cmu}$
are not smooth.   

\sus{Notation} In our calculations below we denote the simple roots by
$\alpha_i, i=0,1,2$; the simple reflection corresponding to $\alpha_i$ is 
denoted by $s_i$; $s_0$ always denotes the affine simple reflection. 

\sus{Type $A_2$}

\sss{} The affine Weyl goup $W_{aff}(A_2)$ of type $A_2$
is described as follows
$$
W_{aff}(A_2)=\{ s_0,\ s_1,\ s_2\ |\ (s_1s_2)^3=1,\ (s_0s_1)^3=1,\ 
(s_0s_2)^3=1\ \} .
$$

There is only one minimal degeneration:

\sss{Minimal singularity of type $a_2$}
This is a singularity of codimension $4$. We have
$$
w=s_1s_2s_1s_0s_1s_2s_1,
\qquad
y=s_1s_2s_1 .
$$
The equivariant multiplicity:
\begin{tiny}
$$
e_{\cla}(\cmu)=\frac{2(3\al_0^2 + 6\al_0\al_1 + 2\al_1^2 + 6\al_0\al_2 
+ 5\al_1\al_2 + 2\al_2^2)}{\al_0 (\al_0 + \al_1) (\al_0 + \al_2) 
(\al_0 + 2 \al_1 + \al_2)(\al_0 + \al_1 + 2 \al_2) 
(\al_0 + 2 \al_1 + 2 \al_2)}
$$
\end{tiny}

\sus{Type $C_2$}

\sss{} The affine Weyl group $W_{aff}(C_2)$ of type $C_2$
is described as follows
$$
W_{aff}(C_2)=\{ s_0,\ s_1,\ s_2\ |\ (s_1s_2)^4=1,\ (s_0s_2)^4=1,\ 
(s_0s_1)^2=1\ \}.
$$

There are two minimal degenerations:

\sss{Minimal singularity of type $c_2$}
$\xymatrix{0\ar@{=>}[r] &1} \leadsto \xymatrix{0\ar@{=>}[r] &0}$ 
This is s singularity of codimension $4$. In this case
$$
w=s_2s_1s_2s_1s_0s_2s_1s_2,
\qquad
y=s_2s_1s_2s_1 .
$$
The equivariant multiplicity:
$$
e_{\cla}(\cmu)=
\frac{8}{\al_0 (\al_0 + 2 \al_2) (\al_0 + 2 \al_1 + 2 \al_2) 
(\al_0 + 2 \al_1 + 4 \al_2)} 
$$

\sss{Quasi-minimal $ac_2$}
$\xymatrix{1\ar@{=>}[r] &1} \leadsto 
\xymatrix{1\ar@{=>}[r] &0}$ 
This is s singularity of codimension $4$. In this case
$$
w=s_2s_1s_2s_1s_0s_2s_0s_1s_2s_0s_2,
\qquad
y=s_2s_1s_2s_1s_0s_2s_0 .
$$
The equivariant multiplicity:
\begin{tiny}
$$
e_{\cla}(\cmu)=\frac{
2(11\al_0^2 + 21\al_0\al_1 + 9\al_1^2 + 43\al_0\al_2 + 
39\al_1\al_2 + 36\al_2^2)}
{\al_0 (\al_0 + \al_1 + \al_2) 
(\al_0 + 2 \al_2) 
(\al_0 + \al_1 + 3 \al_2) 
(2 \al_0 + 3 \al_1 + 4 \al_2) 
(2 \al_0 + 3 \al_1 + 6 \al_2)}
$$
\end{tiny}

\sus{Type $G_2$}\label{G2}

\sss{} The affine Weyl group $W_{aff}(G_2)$ of type $G_2$
is described as follows
$$
W_{aff}(G_2)=\{ s_0,\ s_1,\ s_2\ |\ (s_1s_2)^6=1,\ (s_0s_1)^3=1,\ 
(s_0s_2)^2=1\ \}.
$$

There are three minimal degenerations:

\sss{Minimal singularity of type $g_2$}
$\xymatrix{1\ar@3{>}[r] &0} \leadsto \xymatrix{0\ar@3{>}[r] &0}$.
This is a singularity of codimension 6. We have
$$
w=s_2s_1s_2s_1s_2s_1s_0s_1s_2s_1s_2s_1,
\qquad
y=s_2s_1s_2s_1s_2s_1 .
$$
The equivariant multiplicity:
\begin{tiny}
$$
e_{\cla}(\cmu)=\frac{18}
{\al_0 (\al_0 + \al_1) (\al_0 + \al_1 + 3 \al_2) 
(\al_0 + 3 \al_1 + 3 \al_2) 
(\al_0 + 3 \al_1 + 6 \al_2) (\al_0 + 4 \al_1 + 6 \al_2)}
$$
\end{tiny}

\sss{Quasi-minimal of type $ag_2$}
$\xymatrix{1\ar@3{>}[r] &1} \leadsto \xymatrix{2\ar@3{>}[r] &0}$.
This is a singularity of codimension 4. We have 
$$
w=s_2s_1s_2s_1s_2s_1s_0s_1s_2s_1s_2s_0s_1s_2s_1s_0s_2s_1s_2s_1s_2s_1,
$$
and
$$
y=s_2s_1s_2s_1s_2s_1s_0s_1s_2s_1s_2s_1s_0s_1s_2s_1s_2s_1.
$$
The equivariant multiplicity:
\begin{tiny}
$$
e_{\cla}(\cmu)=\frac{2(27\al_0^2 + 106\al_0 \al_1 + 103\al_1^2 + 
159\al_0\al_2 +309\al_1\al_2 + 216\al_2^2)}
{(\al_0 + \al_1) (\al_0 + 2 \al_1 + 2 \al_2) 
(\al_0 + \al_1 + 3 \al_2) 
 (\al_0 + 2 \al_1 + 4 \al_2) 
(3 \al_0 + 7 \al_1 + 9 \al_2) (3 \al_0 + 7 \al_1 + 12 \al_2)}
$$
\end{tiny}

\sss{Quasi-minimal singularity of type $cg_2$} 
$\xymatrix{0\ar@3{>}[r] &1} 
\leadsto \xymatrix{1\ar@3{>}[r] &0}$.
This is a singularity of codimension 4. We have
$$
w=s_2s_1s_2s_1s_2s_1s_0s_1s_2s_1s_2s_0s_1s_2s_1s_2,
\qquad
y=s_2s_1s_2s_1s_2s_1s_0s_1s_2s_1s_2s_1 .
$$
The equivariant multiplicity:
\begin{equation}\label{cg2mult}
e_{\cla}(\cmu)=\frac{27}
{(\al_0 + \al_1) (\al_0 + \al_1 + 3 \al_2) 
(2 \al_0 + 5 \al_1 + 6 \al_2) (2 \al_0 + 5 \al_1 + 9 \al_2)}
\end{equation}

\sus{Quasi-minimal singularities revisited}

Recall that the singularities $a_2$, $ac_2$, $ag_2$, $c_2$, $cg_2$ have 
codimension 4.
Moreover one observes that the intersection cohomology of singularities
of type $a_2, ac_2, ag_2$ (and, similarly, of $c_2$ and $cg_2$) are the same. 
We conjecture that all these singularities are pairwise smoothly 
non-equivalent.
One verifies that the equivariant multiplicities $e_{\cla}(\cmu)$ are
pairwise distinct (up to linear changes of coordinates with rational 
coefficients). This implies that at least the singularities above are
different as singularities with torus action. 

Similarly, the singularities of types $a_n$ and $ac_n$ have the same 
codimensions and intersection cohomology but we conjecture that these
singularities are not smoothly equivalent. 

\se{Smooth Locus}\label{SmLocSec}

\sus{} Finally we use our results to prove the 
Evens-Mirkovi\' c Theorem 
(Corollary B). 
Let us denote by $\barr\GG^{\sm}_{\cmu}$
the smooth locus of the Schubert variety $\barr\GG_{\cmu}$.

\begin{corollary} For any dominant coweight $\cla\in\CLA^{+}_G$
we have
$$
\barr\GG^{\sm}_{\cmu}=\GG_{\cmu} .
$$
\end{corollary}

\begin{proof} 
It is enough to check that
$\barr\GG_{\cmu}$ is singular along every irreducible
component of the boundary $\barr\GG_{\cmu}-\GG_{\cmu}$.
These irreducible components are precisely Schubert varieties
$\barr\GG_{\cla}$ for all minimal degenerations $\cmu\leadsto\cla$.

Thus we have to check that the minimal degenerations of 
$G(O)$-orbits are singular for all cases in the 
Stembridge's list, see Theorem \ref{Stem}.

In the cases \ref{St1} and \ref{St2} it follows
from the Theorem \ref{LeviTheorem}.

In the cases \ref{sp} and \ref{Stag2} the variety 
$(L^{<0}G\cdot \cla)\cap\barr\GG_{\cmu}$ 
is not rationally smooth by Proposition \ref{ICProposition},
and therefore not smooth.

In the case \ref{Stcg2} the variety 
$(L^{<0}G\cdot \cla)\cap\barr\GG_{\cmu}$ 
is singular by the Kumar's criterion, see formula (\ref{cg2mult}).
\end{proof}

\providecommand{\bysame}{\leavevmode\hbox to3em{\hrulefill}\thinspace}


\begin{thebibliography}{50}

\bibitem{BL}
A.~Beauville and  Y.~ Laszlo, 
\emph{Conformal blocks and generalized theta functions}, 
Comm. Math. Phys. 164 (1994), no. 2, 385--419. 

\bibitem{BD}
A.~ Beilinson and V.~ Drinfeld,
\emph{Quantization of Hitchin's integrable system
and Hecke eigensheaves}, preprint.

\bibitem{Br}
M.~ Brion,
\emph{Equivariant cohomology and equivariant intersection theory}, 
Notes by Alvaro Rittatore. 
NATO Adv. Sci. Inst. Ser. C Math. Phys. Sci.,
514, Representation theories and algebraic geometry 
(Montreal, PQ, 1997), 1--37.

\bibitem{BP}
M.~ Brion and P.~ Polo, 
\emph{Generic singularities of certain Schubert varieties}, 
Math. Z. 231 (1999), no. 2, 301--324.

\bibitem{EM}
S. ~Evens and I.~ Mirkovi\' c, 
\emph{Characteristic cycles for the loop Grassmannian and nilpotent orbits}, 
Duke Math. J. 97 (1999), no. 1, 109--126.

\bibitem{F}
G.~ Faltings,
\emph{A proof for the Verlinde formula},
J. Algebraic Geometry 3 (1994), 347-347.

\bibitem{Fa}
G.~ Faltings,
\emph{Algebraic loop groups and moduli spaces of bundles}, 
J. Eur. Math. Soc. 5 (2003) 1, 41-68.


\bibitem{KP81}
H. ~ Kraft and C.~ Procesi,
\emph{Minimal singularities in ${\rm GL}\sb{n}$}, 
Invent. Math. 62 (1981), no. 3, 503--515. 

\bibitem{KP82}
H. ~ Kraft and C.~ Procesi, 
\emph{On the geometry of conjugacy classes in classical groups}, 
Comment. Math. Helv. 57 (1982), no. 4, 539--602.



\bibitem{Ku2}
S.~ Kumar, \emph{Demazure character formula in arbitrary Kac-Moody setting},
Invent. Math. 89 (1987), no. 2, 395--423.

\bibitem{Ku}
S.~ Kumar, 
\emph{The nil Hecke ring and singularity of Schubert varieties},
Invent. Math. 123 (1996), no. 3, 471--506. 

\bibitem{LS}
Y.~ Laszlo and C.~ Sorger, 
\emph{The line bundles on the moduli of parabolic $G$-bundles 
over curves and their sections}, 
Ann. Sci. \' Ecole Norm. Sup. (4)
30 (1997), no. 4, 499--525. 

\bibitem{Lit}
G.~ Littelmann, \emph{Conracting modules and standard monomial theory for
symmetrizable Kac-Moody algebras}, 
J.~Amer.~Math.~Soc. 11 (1998), no. 3, 551--567. 

\bibitem{L81}
G.~ Lusztig, \emph{Green polynomials and singularities of unipotent 
classes}, 
Adv. in Math. 42 (1981), no. 2, 169--178. 

\bibitem{L83}
G.~ Lusztig, 
\emph{Singularities, character formulas, and a $q$-analog of 
weight multiplicities}, Analysis and topology on singular spaces, 
II, III (Luminy, 1981), 208--229, 
Ast\' erisque, 101-102, Soc. Math. France, Paris, 1983. 


\bibitem{Mat} 
O.~ Mathieu, 
\emph{Formules de caract\` eres pour les alg\` ebres 
de Kac-Moody g\' en\' erales},
Ast\` erisque No. 159-160 (1988), 267 pp. 

\bibitem{MVil}
I.~ Mirkovi\'c and K.~ Vilonen, 
\emph{Perverse sheaves on affine Grassmannians and Langlands duality}, 
Math. Res. Lett. 7 (2000), no. 1, 13--24.

\bibitem{MVyb}
I. ~ Mirkovi\' c and M.~ Vybornov,
\emph{On quiver varieties and affine Grassmannians of type A},
C. R. Acad. Sci. Paris, Ser. I  (2003) 336 (3) 207--212.

\bibitem{MVyb03}
I. ~ Mirkovi\' c and M.~ Vybornov,
\emph{Quiver varieties and Beilinson-Drinfeld 
Grassmannians of type A},
preprint, 2002.


\bibitem{NP}
B.~C.~Ng\^ o and P.~ Polo, 
\emph{R\' esolutions de Demazure affines et formule de 
Casselman-Shalika g\' eom\' etrique}, 
J. Algebraic Geom. 10 (2001), no. 3, 515--547. 

\bibitem{S80}
P.~ Slodowy,  
\emph{Simple singularities and simple algebraic groups}, 
Lecture Notes in Mathematics, 815. Springer, Berlin, 1980. 

\bibitem{St} 
J.~Stembridge, \emph{The partial order of dominant weights},
Adv. Math. {\bf 136} (1998), no. 2, 340-364.

\end{thebibliography}
\end{document}